\documentclass[12pt]{amsart}

\usepackage{amsmath,amsthm,amscd,euscript,amssymb}
\usepackage{graphicx}
\usepackage{epstopdf}
\usepackage{color}

\usepackage[most]{tcolorbox}

\newcommand{\myexample}[2]{
    \begin{tcolorbox}[breakable,colback=black!5!white,colframe=black,title={
    #1}]
        #2
    \end{tcolorbox}
}

\def \c{\mathbb C}
\def \q{\mathbb Q}
\def \z{\mathbb Z}

\def \SL{\hbox{\rm SL}}
\def \id{\hbox{\rm Id}}

\newtheorem{theorem}{Theorem}[section]

\theoremstyle{remark}
\newtheorem{remark}[theorem]{Remark}
\theoremstyle{definition}
\newtheorem{definition}[theorem]{Definition}
\newtheorem{example}[theorem]{Example}
\newtheorem{problem}{Problem}
\newtheorem{conjecture}[problem]{Conjecture}
\newtheorem{question}{Question}

\title[On Hermite's problem]{On Hermite's problem, Jacobi-Perron type algorithms, and Dirichlet groups}

\author{Oleg Karpenkov}

\date{15  January 2021}


\email[Oleg Karpenkov]{karpenk@liverpool.ac.uk}

\begin{document}
\input{epsf}

\begin{abstract}
In 1848  Ch.~Hermite asked if there exists some way to write cubic irrationalities periodically.
A little later in order to approach the problem C.G.J.~Jacobi and O.~Perron
generalized the classical continued fraction algorithm to the three-dimensional case,
this algorithm is called now the Jacobi-Perron algorithm.
This algorithm is known to provide periodicity only for some cubic irrationalities.

In this paper we introduce two new algorithms in the spirit of Jacobi-Perron algorithm:
the heuristic algebraic  periodicity detecting algorithm and the $\sin^2$-algorithm.
The heuristic algebraic  periodicity detecting algorithm is a very fast and efficient algorithm,
its output is periodic for numerous examples of cubic irrationalities,
however its periodicity for cubic irrationalities is not proven.
The $\sin^2$-algorithm is limited to the totally-real cubic case (all the roots of cubic polynomials are real numbers).
In the recent paper~\cite{Karpenkov2021} we proved the periodicity
of the $\sin^2$-algorithm for all cubic totally-real irrationalities.
To our best knowledge this is the first Jacobi-Perron type algorithm for which the cubic periodicity is proven.
The $\sin^2$-algorithm provides the answer to Hermite's problem for the totally real case
(let us mention that the case of cubic algebraic numbers with complex conjugate roots remains open).

We conclude this paper with one important application of Jacobi-Perron type algorithms
to computation of independent elements in the maximal groups of commuting matrices of algebraic irrationalities.
\end{abstract}

\maketitle
\tableofcontents

This paper is dedicated to periodic representations of algebraic numbers.
Recall that a number $\alpha$ is algebraic if
it is a root of some polynomial with integer coefficients.
The smallest degree of integer  polynomial with a root $\alpha$ is called the {\it degree} of the $\alpha$.
It is well known that decimal representations for all rational numbers are eventually periodic or finite,
so the case of algebraic numbers of degree 1 is straightforward.
Let us consider a similar question for algebraic numbers of higher degrees.

\vspace{2mm}

It turns out that the study of this question has a rich history.
Our journey starts
in ancient Greece with the invention of Euclid's algorithm about 300~BC
Euclid's algorithm was originally developed for computing
the greatest common divisor of two integers.
It was two millennia after its invention when the
Euclid's algorithm was used in the study of quadratic irrationals (i.e., algebraic numbers of degree 2).
An important stage here was the introduction of the concept of regular continued fractions
by  J.~Wallis in 1695,
that finally linked Euclid's algorithm to irrational numbers in general and to quadratic irrationalities
in particular.
In 1770 J.-L.~Lagrange proved the periodicity of continued fractions for quadratic irrationalities, closing the question for
the quadratic case  (see Section~\ref{Euclid's algorithm for quadratic irrationalities}).

\vspace{2mm}

For the first time the problem on generalization of Lagrange theorem
on periodicity of continued fractions for quadratic irrationalities to the case of algebraic numbers of degree three
was posed by C.~Hermite in 1849 in a very general settings.
C.~Hermite was wondering if there is a periodic description to cubic irrationalities.
There are many different interpretations of this question that leaded to remarkable theories
in geometry and dynamics of numbers
(see a small survey on various multidimensional generalizations
of ordinary continued fractions in Chapter~23 in~\cite{KarpenkovGCF2013}).

\vspace{2mm}

For this paper we restrict ourselves entirely to the algorithmic approach of the problem that was
initiated by  C.G.J.~Jacobi in 1868 and further developed by O.~Perron in 1907.
They have developed the multidimensional continued fraction algorithm,
now known as the Jacobi-Perron algorithm.
The  Jacobi-Perron algorithm generalizes the Euclidean algorithm and
 provides a sequence of pair of integers similar to a regular continued fractions
provided by the Euclidean algorithm.
The output of the algorithm is periodic for certain cubic numbers,
however  it is believed to be non-periodic for some others.
By that reason Jacobi-Perron algorithm does not provide the complete solution to Hermite's problem,
however it suggests that algorithmic approach might be beneficial to the question.
A similar situation occur with many other Jacobi-Perron type algorithms,
that are neither proven or disproved to be periodic.

\vspace{2mm}

In this paper we introduce two new modifications of Jacobi-Perron
algorithm.
We call the first one the
heuristic algebraic  periodicity detecting algorithm
(or  heuristic APD-algorithm for short) and the second --- the $\sin^2$-algorithm.
The heuristic APD-algorithm is able detect periodicity for all examples in numerous experiments,
it is conjectured to be periodic for algebraic numbers.
The $\sin^2$-algorithm works only in the totally real case (all three roots of the polynomials are real numbers).
For $\sin^2$-algorithm we were able to prove periodicity for triples of cubic conjugate vectors in~\cite{Karpenkov2021}.
(To the best of our knowledge,  this is the first complete proof of periodicity for the Jacobi-Perron type algorithms.)
So the $\sin^2$-algorithm provides the answer to Hermit's problem in the form of Jacobi-Perron type algorithms for
the totally real cubic case.
The non-totally-real case remains open, however we believe that the techniques
of the proof for the $\sin^2$-algorithm can be adapted for that case as well.
Both the heuristic APD-algorithm
and the $\sin^2$-algorithm are discussed in Section~\ref{On periodicity of cubic irrationalities}.

\vspace{2mm}

We say just a few words regarding higher dimensional case in Section~\ref{Situation in degree greater than 3},
which currently remains open.

\vspace{2mm}

Further we address an important application of Jacobi-Perron type algorithms.
It turns out that such algorithms provide a simple way to write
independent (with respect to the matrix multiplication operation) commuting pairs of matrices for the corresponding Dirichlet group.
Recall  that the groups of commuting matrices (so called Dirichlet groups)
are described by a mysterious Dirichlet's unit theorem (we formulate and discuss it later in Subsection~\ref{Dirichlet's unit theorem}),
whose complete understanding will probably bing light on the periodicity of generalized Euclidean algorithms.
Classical proofs of Dirichlet's unit theorem provide huge estimates on the coefficients of the generators of the Dirichlet groups.
The brute force algorithms provided by this theorem are very slow and seems to  have no practical value.
We discuss a simple and fast approach to the problem in the last two chapters of the paper.

\vspace{2mm}

\section{Euclid's algorithm for quadratic irrationalities}
\label{Euclid's algorithm for quadratic irrationalities}

As we have already mentioned
the periodicity of  quadratic irrationalities is closely related to Euclid's algorithm.
Recall that the original algorithm computes the greatest common divisor of
two integer numbers.
Let us first write down a slightly extended form such that it can be applied to arbitrary
numbers (not necessarily integers).

\vspace{2mm}


\myexample{Extended Euclid's algorithm}
{
{\noindent{\bf \underline{Input}:}}
We start with  a pair o real numbers $(p, q)=(p_0,q_0)$ assuming $q_0>0$.

\vspace{1mm}

{\noindent {\bf \underline{Step of the algorithm}:}}
Assume that we have found two real numbers $(p_i, q_i)$ with $q_i\ge 0$.
Then the next step is:
$$
(p_i,q_i)\mapsto
(p_{i+1},q_{i+1})=(q_i, p_i-\lfloor p_i/q_i\rfloor q_i)
$$
Here we call the value $a_i=\lfloor p_i/q_i\rfloor$ as the {\it $i$-th element} of the algorithm.

\vspace{1mm}

{\noindent{\bf \underline{Termination of the algorithm}:}}
In case if we have arrived to the pair  $(a_i, b_i=0)$
we do not proceed further. Here the algorithm terminates.
}

\begin{remark}
It is interesting to note that in the case of a pair of integers $(p,q)$ with $q>0$
we have the classical Euclid's algorithm.
Here the algorithm terminates in the finite number of steps and on the last step we get
$(\gcd(p,q),0)$ where $\gcd(p,q)$ is the greater common divisor of $p$ and $q$.
\end{remark}

\begin{remark}
The algorithm terminates if $p/q$ is a rational number,
and it does not terminate otherwise.
\end{remark}

\begin{example}
Let us apply the algorithm to the pair $(21,15)$. We have:
$$
(21,15) \mapsto (15,6)\mapsto(6,3)\mapsto (3,0).
$$
The output of the algorithm is as follows:
$$a_1=1, \quad a_2=2, \quad \hbox{and} \quad  a_3=2.
$$
Note that
$$
\gcd(21,15)=3 \quad \hbox{and} \quad
\frac{21}{15}=1+\frac{1}{2+1/2}.
$$
\end{example}

Now let us focus  on the case of pairs $(\alpha,1)$ where $\alpha$ is any real number.
In this case the extended Euclid's algorithm generates a remarkable sequence of numbers $a_i$.
If $\alpha$ is a rational number, then  the algorithm terminates (on the $n$-th step for some integer $n$)
and the output sequence $(a_i)$ satisfies the following identity.
$$
\alpha=a_0+\frac{1}{\displaystyle a_1+\frac{1}{\displaystyle
a_2+\frac{1}{\displaystyle\ddots+\frac{\displaystyle 1}{a_n}}}}
$$
The expression from the right hand side is called a {\it regular continued fractions} for $\alpha$
and denoted by  $[a_0;a_1:\cdots :a_n]$.
(Notice here that the term continued fraction was introduced by J.~Wallis in 1695.)

\vspace{2mm}

The above identity for rational $\alpha$ is extended to the case
of irrational $\alpha$ by the  following limit
$$
\lim\limits_{k\to \infty} [a_0;a_1:\cdots: a_k],
$$
which we call the {\it regular continued fraction} for $\alpha$ and denote by $[a_0;a_1:\cdots]$.
Let us just notice that this limit always exists and
distinct sequences converge to distinct irrational numbers.
(For the details of the classical theory of continued fractions we refer, e.g., to~\cite{Khinchin1961}.)

\vspace{2mm}

We are finally arriving to a very non-trivial theorem on periodicity of continued fractions for quadratic irrationalities.
This theory introduced almost a century later in 1770 by J.-L.~Lagrange (see in Chapter 34 of~\cite{Lagrange1770}).

\begin{theorem}{\bf (J.-L.~Lagrange.)}\label{lagr}
A regular continued fraction of $\alpha$ is periodic
if and only if $\alpha$ is a quadratic
irrationality $($i.e., $\frac{\displaystyle
a+b\sqrt{c}}{\displaystyle d}$ for some integers $a$, $b$, $c$, and
$d$, where $b\ne 0$, $c> 1$, $d>0$, and $c$ is square-free$)$.
\end{theorem}

This theorem gives a complete answer to the question on periodic representations for quadratic irrationalities.

\begin{example}\label{2d-periodic-ex}
Let us apply the Extended Euclid's algorithm to $(2\sqrt{5},1)$.
We have
$$
(2\sqrt{5},1) \mapsto c_1(1+\sqrt{5}/2,1) \mapsto c_2(4+2\sqrt{5},1) \mapsto c_3(1+\sqrt{5}/2,1)\mapsto\ldots
$$
where
$$
c_1=2\sqrt{5} - 4, \quad c_2=9 - 4\sqrt{5}, \quad c_3=34\sqrt{5} - 76, \quad \ldots
$$
Note that the vectors obtained on the first and on the third step are proportional.
Hence the output of the Euclidean algorithm is periodic with 1 element in period and 2 elements in pre-period.
Here we have
$$
a_1=4, \quad a_{2k}=2, \quad \hbox{and} \quad  a_{2k+1}=8
$$
for all integer  $k>1$.
Note that
$$
2\sqrt{5}=[4;2:8:2:8:2:8:\ldots].
$$
\end{example}

\section{On periodicity of cubic irrationalities}
\label{On periodicity of cubic irrationalities}

The problem of detecting periodicity for cubic irrationalities was
posed by Ch.~Hermite in 1848 (see e.g.~\cite{Picard1901}, \cite{Hermite1850}), where
he was asking wether
{\it there exists some way to write cubic irrationalities periodically?}
In this section we discuss the Jacobi-Perron algorithmic approach and recent advances in it.

\subsection{Jacobi-Perron algorithm}

Jacobi-Perron algorithm is one of the possible ways to generalize extended Euclid's algorithm to higher dimensions,
it was proposed by
C.G.J.~Jacobi in 1868 in~\cite{Jacobi1868} and further developed by
O.~Perron in 1907, see~\cite{Perron1907}.
The algorithm is as follows.

\vspace{2mm}

%
%

\myexample{Jacobi-Perron algorithm}
{

{\noindent{\bf \underline{Input}:}}
We starts with triples of real numbers $(x,y,z)$.

\vspace{1mm}

{\noindent {\bf \underline{Step of the algorithm}:}}
In the previous step we have constructed $(x_i, y_i,z_i)$.
Then we proceed with the following iteration:
$$
(x_i,y_i,z_i) \mapsto (x_{i+1},y_{i+1},z_{i+1})=\Big(y_i,z_i-\Big\lfloor\frac{z_i}{y_i}\Big\rfloor y_i, x_i-\Big\lfloor\frac{x_i}{y_i}\Big\rfloor y\Big).
$$
Here the {\it $i$-th element} of the corresponding multidimensional continued fraction is set to be the
pair of integers
$$
\Big(\Big\lfloor\frac{z_i}{y_i}\Big\rfloor,\Big\lfloor\frac{x_i}{y_i}\Big\rfloor\Big).
$$
\vspace{1mm}

{\noindent{\bf \underline{Termination of the algorithm}:}}
If we have arrived to the triple  $(x_i, y_i,z_i)$ where $y_i=0$, then
the algorithm terminates.
}


\vspace{2mm}

\begin{remark}{\bf (How to generate cubic vectors.)}\label{cubic vectors}
As we have seen the input data for the Jacobi-Perron algorithm is a triple of numbers.
Let  us discuss how to write a cubic vector starting from a cubic number $\alpha$.
For the first two coordinates of this vector we take 1 and $\alpha$.
Now it remains to find out how to pick the last coordinate of this vector.
There is a natural answer to this question.
Consider an arbitrary polynomial $q$ of degree 2 with integer coefficients and
let us take the vector
$$
(1,\alpha,q(\alpha)).
$$
The simplest here would be $(1,\alpha,\alpha^2)$, here we take the polynomial $q(x)=x^2$.

\vspace{1mm}

In general, one can pick three numbers in $\q(\alpha)$ that form a basis
of the linear space $\q(\alpha)$ over $\q$.
Of course, the choice of the basis of  $\q(\alpha)$ will result in different outputs of the Jacobi-Perron algorithm.
The problem of {\it describing all possible periods for continued fraction algorithms for different vectors of $\q(\alpha)$}
 remains open {\it for every single $\alpha$}.
In particular it is not known what is the set of available periods
 for the classical case of regular continued fractions of quadratic irrationalities
 in $\q(\sqrt{2})$, $\q(\sqrt{3})$, $\q(\sqrt{5})$, etc.
\end{remark}

Let us continue with the following example.

\begin{example}\label{ex-ok}
Let $\xi$ be a real root of the polynomial $x^3+2x^2+x+4$,
namely
$$
\xi=
-\frac{(53 + 6\sqrt{78})^{1/3}}{3} - \frac{1}{3(53 + 6\sqrt{78})^{1/3}} - \frac{2}{3}.
$$
Now consider a vector
$$
(1,\xi,\xi^2+\xi).
$$

Then the Jacobi-Perron algorithm will generate the following periodic output.
\begin{center}
\begin{tabular}{|c||c|c|c|c|c|c||c|c||c|c||c|c||c|c|c|c|c|c|c|c|c|c|c}
\hline
&1 &2 &3 &4 &5 &6 &$2k+1$ &$2k+2$
\\
\hline
\hline
$\lfloor x/y\rfloor$ & -1 &1 &1 &1 &2 &6 &3 &7
\\
\hline
$\lfloor z/y\rfloor$& -2 &0 &0 &0 &2 &4 &1 &1
\\
\hline
\end{tabular}
\end{center}
(Here $k\ge 3$.) After the first 6 steps of the algorithm the sequence starts to be periodic with period 2.
\end{example}

The question of periodicity for the Jacobi-Perron algorithm is known in folklore
as Jacobi's Last Theorem.

\begin{problem}{\bf (Jacobi's Last Theorem.)}
Let $K$ be a totally real cubic number field.
Consider arbitrary elements $y$ and  $z$ of $K$
satisfying $0<y,z<1$ such that $1$, $y$, and $z$ are independent elements over $\q$.
Is it true that the Jacobi-Perron algorithm generates an eventually periodic continued fraction with
starting data $v=(1,y,z)$?
\end{problem}

The answer to the question of Jacobi's Last Theorem
seems to be negative. Let us consider another example to see this.

\begin{example}\label{not-ok}
Let us consider the vector
$$
v=(1,\sqrt[3]{4},\sqrt[3]{16}) .
$$
Numerical computations suggest that the output of the Jacobi-Perron algorithm for this vector
is not eventually periodic. Here we show the output elements for the first several steps of the algorithm
(for further numerical computations and discussions we refer to~\cite{Elsner1967}).

\begin{center}
\begin{tabular}{|c||c|c|c|c|c|c|c|c|c|c|c|c|c|c|c|c|c|c|c|c|c|c|c}
\hline
&1 &2 &3 &4 &5 &6 &7 &8 &9 &10 &11 &12 &$\ldots$&94&$\ldots$
\\
\hline
\hline
$\lfloor x/y\rfloor$ & 0 &1 &13 &1 &6 &1 &1 &3 &2 &3 &4 &1 &$\ldots$&476&$\ldots$
\\
\hline
$\lfloor z/y\rfloor$& 1 &1 &9 &1 &2 &0 &0 &2 &0 &1 &1 &1 &$\ldots$&388&$\ldots$
\\
\hline
\end{tabular}
\end{center}

\end{example}

\subsection{A few words about Gauss-Kuzmine statistics}\label{A few words about Gauss-Kuzmine statistics}
Let us make a small informal discussion of the last example.
The sequence of the last example seems to be non-periodic.
One could notice that the numbers of the elements are relatively small,
we have a few bumps only. In fact, this is a rather predictible behaviour for non-periodic
sequences.
Let us consider the first 30 digits for the regular continued fraction for $\pi$. We have
$$
\pi=[3: 7; 15; 1; 292; 1; 1; 1; 2; 1; 3; 1; 14; 2; 1; 1; 2; 2; 2; 2; 1; 84; 2; 1; 1; 15; 3; 13; 1; 4; 2; \ldots]
$$
As we see, the most frequent element is 1; the next frequent element is 2. etc.
This phenomenon is described by the Gauss-Kuzmine theorem  stating that the frequency
of an element $k$ is
$$
\frac{1}{\ln(2)}\ln\bigg(1+\frac{1}{k(k+1)}\bigg).
$$
For the first time it was proved by R.O.~Kuzmin in 1928 in~\cite{Kuzmin1928} (see also in~\cite{Kuzmin1932}).
It is interesting to notice that the Gauss-Kusmine statistics has a projective nature, it can be written
in terms of cross-rations:
$$
\frac{1}{\ln(2)}\ln\bigg(1+\frac{1}{k(k+1)}\bigg)=\frac{\ln[-1,0,k,k+1]}{\ln[-1,0,1,\infty]}.
$$
It remains to say that
{\it the analogs of Gauss-Kuzmine theorem
for Jacobi-Perron type algorithms in higher dimensions are not known},
however we might expect a similar behaviour for the elements in higher dimensions as well.

For further discussions and first successful generalisation of Gauss-Kuzmine theorem to higher dimensional case
we refer to~\cite{Kontsevich1999} and~\cite{Karpenkov2007}, see also Chapter 19 in~\cite{KarpenkovGCF2013}.

\subsection{Heuristic algebraic periodicity detecting algorithm}

Computations of  L.~Elsner and H.~Hasse~\cite{Elsner1967} suggest
that the output of the Jacobi-Perron algorithm for the cubic vector $(1,\sqrt[3]{4},\sqrt[3]{16})$  of Example~\ref{not-ok}
is non-periodic.
However  {\it the proof of this fact is missing},
there is a strong believe that the sequence is indeed not periodic.

\vspace{2mm}

Let us informally say a few words on the reason for the Jacobi-Perron algorithm potentially
to be non-periodic for cubic vectors.
In fact, any cubic vector has a pair of algebraically conjugate vectors that are
completely defined by the original vector.
The pairs of these three vectors generate an arrangement of three planes  with an action of
the corresponding Dirichlet groups that we will discuss in Section~\ref{Dirichlet groups}.
In some sense  the choice of the elements in the classical Jacobi-Perron algorithm
is blind to the action of the corresponding Dirichlet group, it follows
more the Euclidean distances to nearest integers. The last seems to be not appropriate
for cubic vectors.

\vspace{2mm}

Let us introduce an important ternary form related to triples of vectors.
It will be further used for triples of cubic conjugate vectors.

\begin{definition}
Consider three vectors
$$
u=(u_1,u_2,u_3), \quad
v=(v_1,v_2,v_3), \quad \hbox{and} \quad
w=(w_1,w_2,w_3)
$$
in $\c^3$.
The following ternary form
$$
\det
\left(
\begin{array}{ccc}
x& y& z\\
v_1& v_2& v_3\\
w_1& w_2& w_3\\
\end{array}
\right)
\cdot
\det
\left(
\begin{array}{ccc}
u_1& u_2& u_3\\
x&y&z\\
w_1& w_2& w_3\\
\end{array}
\right)
\cdot
\det
\left(
\begin{array}{ccc}
u_1& u_2& u_3\\
v_1& v_2& v_3\\
x& y& z\\
\end{array}
\right)
$$
in three variables $x$, $y$, and $z$
is called the {\it Markov-Davenport characteristics}
of $(x,y,z)$ with respect to the vectors $u,v,w$.
Denote it by $\chi_{u,v,w}.$
\end{definition}

\begin{remark}
The Markov-Davenport characteristics was first studied in the context of their minima
in series of works~\cite{Davenport1938, Davenport1938a,Davenport1939,Davenport1943}
by H.~Davenport in the middle of the XX century.
These minima generalises two-dimensional
Markov minima introduced by A.~Markov in
1879 in~\cite{Markoff1879}
(for more details see a very nice book~\cite{Aigner2013} by M.~Aigner and also the paper~\cite{KV2020}).

\end{remark}

\begin{remark}\label{c3}
Note that the set of zeroes for the Markov-Davenport characteristics
is the union of all invariant planes in $\c^3$.
The  Markov-Davenport characteristics shows
{\it how close is the point from the union of planes spanned
by pairs of vectors $(u,v)$, $(v,w)$, and $(w,u)$}.
\end{remark}

Now we would like to introduce a modification of the Jacobi-Perron algorithm
that will be aiming to minimize Markov-Davenport characteristics (rather than the Euclidean distance to the nearest integer vector).
We would like to call this algorithm the
{\it heuristic algebraic  periodicity detecting algorithm} or the {\it heuristic APD-algorithm}, for short.

\vspace{2mm}

%
%

\myexample{Heuristic APD-algorithm}
{

{\noindent{\bf \underline{Input}:}}
 One starts with triples of real vectors $(\xi,\nu,\mu)$ where the last coordinate
of $\xi=(x_0,y_0,z_0)$ is positive (i.e., $z_0>0$),

\vspace{1mm}

{\noindent
{\bf \underline{Step~$0$}:}} First of all let us make all the coordinates of $\xi$ to be positive by applying
the following integer lattice preserving transformation:
$$
T_0: (x,y,z) \mapsto
\Big(x-\Big\lfloor\frac{x}{z}\Big\rfloor z, y-\Big\lfloor\frac{y}{z}\Big\rfloor z,z\Big),
$$
namely we consider
$$
(\xi_1,\nu_1,\mu_1)= \big(T_0(\xi),T_0(\nu),T_0(\mu)\big).
$$

{\noindent
{\bf \underline{Step~$i$ for $i\ge1$}}:} In the previous step we have constructed $(\xi_i,\nu_i,\mu_i)$
with positive coordinates of $\xi_i=(x_i,y_i,z_i)$.
In additional from Step~2 on we have $y_i>x_i$ and $y_i>z_i$.

\vspace{1mm}

\begin{itemize}
\item {\bf Stage 1: Determination of the element of continued fraction:} The new element of the heuristic continued fraction  $(a_i,b_i)$ is defined from the following four conditions:

--- $0\le  a_i \le \lfloor x_i/z_i\rfloor$;

--- $0\le  b_i \le \lfloor y_i/z_i\rfloor$;

--- $(a_i,b_i)\ne (0,0)$ (with the only one exception for Step~1: we choose $(0,0)$ if  $\lfloor x_1/z_1\rfloor<1$ and $\lfloor y_1/z_1\rfloor<1$);

--- the triple  $(a_i,b_i,1)$ provides the minimal possible value for the absolute value of the Markov-Davenport characteristic
$\chi_{\xi,\nu,\mu}$
among the vectors $(a,b,1)$ whose first two coordinates satisfy the first three conditions.

\vspace{1mm}

\item {\bf Stage 2: Iteration step:}
$$
T_i:
(x,y,z) \mapsto (y-b_i z,z,x-a_iz).
$$
Here we construct
$$
(\xi_{i+1},\nu_{i+1},\mu_{i+1})=\big(T_i(\xi_i),T_i(\nu_i),T_i(\mu_i)\big)
$$
\end{itemize}

\vspace{1mm}

{\noindent{\bf\underline{Termination of the algorithm}.}}
In the case if the last coordinate of $\xi_i$ is zero (i.e. $z_i=0$)
we do proceed further and the algorithm terminates.
}


\begin{remark}
As one can see the iteration step of the heuristic APD-algorithm
reminds the Jacobi-Perron algorithm.
Here the main difference between the algorithms is as follows.
The Jacobi-Perron algorithm takes maximal possible values for $a_i$ and $b_i$,
which would the best to approximate coordinatewise.
However as Example~\ref{not-ok} shows, the coordinatewise approximation is not a good approximation
with respect to Markov-Davenport characteristic.
In the heuristic APD-algorithm we are aiming to minimize  Markov-Davenport characteristic,
which is due to Remark~\ref{c3} (informally speaking) provides a simultaneous approximation.
\end{remark}

\begin{remark}
The idea of the heuristic APD-algorithm appeared during the study of
Klein polyhedra by the author.
Notice that Klein polyhedra were introduced in 1895 by F.~Klein in~\cite{Klein1895,Klein1896} (roughly at the time
when the Jacobi-Perron algorithm appeared for the first time).
Theory of Klein polyhedra represents the combinatorial periodicity of algebraic cones.
Klein polyhedra are known to be doubly periodic for the case of totally real cubic numbers.
There is no much known about the link between periodicity of Klein polyhedra and Jacobi-Perron type algorithms,
however they are both related to Dirichlet groups discussed briefly below.
(For further details on Klein polyhedra we refer to~\cite{Arnold2002,Arnold1998,KarpenkovGCF2013}.)

\end{remark}

Heuristic APD-algorithm is designed to work with
triples of cubic conjugate vectors.
For simplicity we define conjugate vectors using the following property.

\begin{definition}\label{conjugate}
Let $M$ be a matrix with integer elements, and let the characteristic polynomial of $M$
be irreducible over $\q$.
Then the triples of distinct eigenvectors of $M$ are said to be {\it cubic conjugate vectors}.
\end{definition}

\begin{remark}{\bf (How to generate triples of cubic conjugate vectors from a single cubic number.)}\label{conj-remark}
Let $\alpha$ be a cubic number and let $p$ be any polynomial with integer coefficients of degree $3$ having $\alpha$ as a root.
(We assume that $p$ is irreducible over $\q$.)
Similarly to the case of single cubic vectors
(see Remark~\ref{cubic vectors})
we can construct triples of conjugate cubic vectors.
In order to do this we
additionally pick an arbitrary degree 2 polynomial $q$ with integer coefficients.
Now a triple of conjugate vectors can be naturally derived from $(\alpha,p,q)$.
Namely, let $\alpha$, $\beta$, and $\gamma$ be distinct roots of $p$.
Then we set
$$
\begin{array}{l}
\xi=\big(1,\alpha,q(\alpha)\big);\\
\nu=\big(1,\beta,q(\beta)\big);\\
\mu=\big(1,\gamma,q(\gamma)\big).\\
\end{array}
$$
It turns out that these vectors are eigenvectors of some integer matrix, and hence they are cubic conjugate vectors.
\end{remark}

\begin{example}\label{not-ok-2}
Let us now consider the cubic vector of Example~\ref{not-ok} for which we have failed to detect the periodicity
of the Jacobi-Perron algorithm output:
$$
\xi=(1,\sqrt[3]{4},\sqrt[3]{16}) .
$$
Note that  $\sqrt[3]{4}$ is a root of $x^3-4$.
Note also that
$$
\sqrt[3]{16}=(\sqrt[3]{4})^2
$$
Let $\beta$ and $\gamma$ be two other complex roots of $x^3-4$.
Consider two vectors:
$$
\nu=(1,\beta,\beta^2) \quad \hbox{and} \quad \mu =(1,\gamma,\gamma^2).
$$

Then the output of the heuristic APD-algorithm for the triple of vectors $(\xi,\nu,\mu)$ is as follows.

\begin{center}
\begin{tabular}{|c||c||c|c|c|c|c|c||c|c|c|c|}
\hline
&0 &1 &2 &3 &4 &5 &6 &$4k+3$ &$4k+4$ &$4k+5$ &$4k+6$
\\
\hline
\hline
$a_1$ & 0 &0 &1 &0 &0 &0 &1 &1 &1 &0 &0
\\
\hline
$b_1$ & 0 & 0 &2 &1 &1 &1 &5 &0 &1 &1 &6
\\
\hline
\end{tabular}
\end{center}
(Here $k \ge 1$.)

\vspace{2mm}

Note that Step~0 does not change the triple.
At Step~1 we have a situation when the last element is the greatest, so we have the exception $(a_1,b_1)=(0,0)$.
From Step~2 on the greatest coordinate will be always the second one, so we will never have $(0,0)$ as an element
of current step.
After 6 steps of a pre-period we have a periodic sequence with period 4.
\end{example}

Let us continue with the following conjecture.

\begin{conjecture}
The heuristic APD-algorithm is periodic for all triples of cubic conjugate vectors
of Definition~\ref{conjugate}.
\end{conjecture}

\begin{remark}
The conjecture can be considered separately for both totally-real and complex cases.
This conjecture has a straightforward generalization to the case of $d$-tuples of conjugate
algebraic vectors  of degree $d$ (see some further discussions in Section~\ref{Situation in degree greater than 3}).
\end{remark}

\subsection{Technical remarks on  the heuristic APD-algorithm}

While practically approaching particular examples of quadratic numbers,
on can use the following two approaches.
The first one is more symbolical: here
we represent quadratic numbers using quadratic root expression and do all the computations with them
(e.g., see Example~\ref{2d-periodic-ex}).
The second approach is to write sufficiently precise rational decimal approximations of quadratic numbers and find the period
using their continued fractions.

Similar to quadratic case, cubic and quartic cases have a symbolic approach due to Cardano formula
and  Ferrari's method.
It is also simple to detect periodicity having a sufficiently good approximation of cubic vectors.
For instance, in order to find periodic representation for the vector of Example~2,
it is enough to know only the first 5 decimal digits of the coordinates of this vector.
In general, the smaller the periods and  the pre-period are and the smaller the elements of the pre-period and the period are
the smaller the rate of approximation is sufficient for constructing the period and the pre-period.

\begin{remark} (On Gauss-Kuzmine statistics.)
Practically the algorithm is very fast if the obtained elements of continued fractions
are small, and it is starting to be slower with the growth of the elements.
Here we should mention that in practice large elements occur very seldom,
so on average algorithm works fast.
As in the case of the Jacobi-Perron algorithm (see Subsection~\ref{A few words about Gauss-Kuzmine statistics})
the {\it distribution of frequencies is unknown for the case of the  heuristic APD-algorithm}.
\end{remark}

\subsection{A remark on the periodicity of  the $\sin^2$-algorithm in the cubic  totally real case}

It seems that heuristic APD-algorithm detects periodicity of cubic vectors.
In addition it works very fast, so heuristically it solves effectively the problem of finding the pre-periods and the periods
for cubic vectors.
Thus practically it can be effectively used for the computations of the independent elements in Dirichlet groups
(or units in orders of algebraic fields), we discuss this later in Section~\ref{Dirichlet groups}.
Currently, the main concern regarding this  algorithm is that proof for cubic periodicity is missing for it.

\vspace{2mm}

Recently we have developed an algorithm similar to the heuristic APD-algorithm
that is designed for the totally real case (i.e., when the corresponding cubic extension of rational numbers is
embedable to the real line)
and proved its periodicity (see~\cite{Karpenkov2021}).
As far as we are aware of, this is the first complete proof of periodicity for the Jacobi-Perron type algorithms.

\vspace{2mm}

Let us outline the idea of the algorithm.  Given three real vectors $(\xi,\nu,\mu)$,
in each step of  the algorithm we aim to minimize the $\sin^2$
of the angle between the planes spanned by pairs of vectors $(\xi,\nu)$ and $(\xi,\mu)$.
One can say that the value of $\sin^2$ here replaces the Markov-Davenport characteristics
in the heuristic APD-algorithm.

\vspace{2mm}

%
%

\myexample{$\sin^2$-algorithm}
{

{\noindent{\bf \underline{Input}:}}
 We are given by three vectors $\xi$, $\nu$, $\mu$ such that

 \vspace{1mm}

 --- $\xi$ has positive coordinates $(x,y,z)$ satisfying $x>y>z>0$;

 --- all three coordinates for both $\nu$ and $\mu$ are neither simultaneously positive nor simultaneously negative.

\vspace{2mm}

{\noindent {\bf \underline{Step of the algorithm}:}}
Let us apply the following linear transformation
$$
(\xi_i, \nu_i, \mu_i) \to (\Phi_i(\xi_i),\Phi_i(\nu_i),\Phi_i(\mu_i))
$$
with
$$
\Phi_i=T_iM_i.
$$
Here $M_i$ is taken to be the minimizer of the value
of $\sin^2$ for the angle between the plane spanned by
$M_i(\xi_i)$  and $M_i(\nu_i)$ and
the plane spanned by $M_i(\xi_i)$ and $M_i(\mu_i)$.
The minimization is done among all the transformations
$$
N_{\alpha,\beta,\gamma}: (x,y,z)\mapsto \big(x-\alpha z-\gamma(y-\beta z), y-\beta z, z)
$$
satisfying
$$
0\le \alpha\le \Big\lfloor \frac{x_i}{z_i}\Big\rfloor, \quad
0\le \beta \le \Big\lfloor \frac{y_i}{z_i}\Big\rfloor,  \quad \hbox{and} \quad
0\le \gamma \le \Big\lfloor \frac{x_i/z_i-\alpha}{y_i/z_i-\beta}\Big\rfloor,
$$
and the transformation
$$
N_0= (x,y,z)\mapsto  \big(x - y, y, z-(x-y)\big),
$$
which is considered only in case if
$$
z_i>x_i-y_i>0.
$$

After $M_i$ is constructed we set  $T_i$ as a transposition of the basis vectors that puts the coordinates
of $M_i(\xi)$ in the decreasing order.

\vspace{2mm}

At each step the algorithm returns $\Phi_i$ is an output.

\vspace{2mm}

{\noindent{\bf \underline{Termination of the algorithm}:}}
In the case if the last coordinate of $\xi_i$ is zero (i.e. $z_i=0$)
we do proceed further and the algorithm terminates.

}

\begin{remark}
Note that the triples in input of the algorithm have some initial conditions.
There is a rather simple way to change the coordinates of an arbitrary triple of vectors
 (by an integer lattice preserving linear transformation)
such that in the new basis this triple fulfills the input conditions.
We omit the technical details here.
For further details and the proof of periodicity for cubic vectors
we refer the interested reader to~\cite{Karpenkov2021}.
\end{remark}

\begin{remark}
Finally we would like to refer to several research papers on cubic periodicity in some other settings.
Cubic periodicity was also studied for the cases of the following generalized continued fractions:
for Klein polyhedra~\cite{Lachaud1993,German2008},
Minkovski-Voronoi polyhedra~\cite{Voronoi1952,Minkowski1967,Bullig1940},
triangle sequences~\cite{Dasaratha2012}, and
ternary continued fractions (or bifurcating continued fractions)~\cite{Murru2015}.
\end{remark}

\section{Situation in degree greater than 3}\label{Situation in degree greater than 3}

This section is rather small as it is {\it almost nothing known in the cases of degree greater than 3}.
Currently the main source of ideas that are applied in the  higher degree case
are coming from the  study of cubic vectors.

\vspace{2mm}

We should notice that {\it there is a straightforward generalization of the heuristic APD-algorithm, which is likely to detect periodicity for irrationalities of degree $d>3$}.
Let us briefly formulate it.

\vspace{2mm}
Consider an irreducible polynomial of degree $d$ let $\xi_1,\ldots \xi_d$
be the set of its solutions.
Let $q_1,\ldots,q_d$ be a basis of polynomials of degree less than $d$ with rational coefficients.

consider vectors
$$
(q_1(\xi_i),q_2(\xi_i),\ldots, q_d(\xi_i)) \quad \hbox{ for $i=1,\ldots, d$.}
$$
The Markov-Davenport characteristic for these vectors is now written as a product of $d$ matrices of size $d\times d$:
$$
\xi(x_1,\ldots, x_d)=\prod\limits_{k=1}^{d} M_k
$$
where $M_k(\xi)$ is the obtained from matrix
$$
\left(
\begin{array}{cccc}
q_1(\xi_1)&q_1(\xi_2)&\ldots&q_1(\xi_d)\\
q_2(\xi_1)&q_2(\xi_2)&\ldots&q_2(\xi_d)\\
\vdots&\vdots&\ddots&\vdots\\
q_d(\xi_1)&q_d(\xi_2)&\ldots&q_d(\xi_d)\\
\end{array}
\right)
$$
by replacing its $k$-th column by the column of variables $(x_1,\ldots, x_d)$.
It is interesting to note that after a multiplication by some constant
the coefficients of the Markov-Davenport characteristic are all integers, see, e.g., in Chapter 21.4 in~\cite{KarpenkovGCF2013}).

\vspace{2mm}

The multidimensional heurostic APD-algorithm will be as follows.

\vspace{2mm}

%
%

\myexample{Multidimensional heuristic APD-algorithm}
{

{\noindent{\bf \underline{Input}:}}
 One starts with triples of conjugate real vectors $(\xi_{1,0},\xi_{2,0},\ldots, \xi_{d,0})$
 generated as above, where the last coordinate
of $\xi_{1,0}=(x_{1,1,0},x_{1,2,0},\ldots, x_{1,d,0})$ is positive (i.e., $x_{1,d,0}>0$),

\vspace{1mm}

{\noindent
{\bf \underline{Step~$0$}:}} First of all let us make all the coordinates of $\xi_{1,0}$ to be positive by applying
the following integer lattice preserving transformation:
$$
T_0: (x_1,x_2,\ldots, x_d) \mapsto
\Big(x_1-\Big\lfloor\frac{x_1}{x_d}\Big\rfloor x_d, x_2-\Big\lfloor\frac{x_2}{x_d}\Big\rfloor x_d,
\ldots, x_2-\Big\lfloor\frac{x_{d-1}}{x_d}\Big\rfloor x_d, x_d\Big),
$$
namely we consider
$$
(\xi_{1,0},\xi_{2,0},\ldots, \xi_{d,0})= \big(T_0(\xi_{1,0}),T_0(\xi_{2,0}),\ldots, T_0(\xi_{d,0})\big).
$$

{\noindent
{\bf \underline{Step~$i$ for $i\ge1$}}:} In the previous step we have constructed
the $d$-tuple
$$
(\xi_{1,i},\xi_{2,i},\ldots, \xi_{d,i})
$$
with positive coordinates of $\xi_{1,i}=(x_{1,1,i},x_{1,2,i},\ldots, x_{1,d,i})$.
In additional from Step~2 on we have $x_{1,d-1,i}>x_{1,k,i}$ for all $k<d-1$ and $k=d$.

\vspace{1mm}

\begin{itemize}
\item {\bf Stage 1: Determination of the element of continued fraction:} The new element of the heuristic continued fraction  $(a_{1,i},a_{2,i},\ldots, a_{d-1,i})$ is defined from the following conditions:

--- $0\le  a_{k,i} \le \lfloor x_{1,k,i}/x_{1,d,i}\rfloor$ for $k=1,\ldots, d-1$;

--- $(a_{1,i},a_{2,i},\ldots, a_{d-1,i})\ne (0,0,\ldots, 0)$ (with the only one exception for Step~1: we choose $(0,0)$ if  $\lfloor x_{1,k,i}/x_{1,d,i}<1$ for $k=1,\ldots, d-1$);

--- the $d$-tuple  $(a_{1,i},a_{2,i},\ldots, a_{d-1,i},1)$ provides the minimal possible value for the absolute value of the Markov-Davenport characteristic
among the triples  $(a_1,\ldots,a_d,1)$ whose first two coordinates satisfy the first two conditions.

\vspace{1mm}

\item {\bf Stage 2: Iteration step:}
$$
T_i:
(x_1,x_2,\ldots, x_d) \mapsto (x_2-a_{2,i} x_d,\ldots,x_{d-1}-a_{d-1,i} x_d,x_d,x_1-a_{1,i}x_d).
$$
Here we construct
$$
(\xi_{1,i+1},\xi_{2,i+1},\ldots, \xi_{d,i+1})=\big(T_i(\xi_{1,i}),T_i(\xi_{2,i}),\ldots, T_i(\xi_{d,i})\big)
$$
\end{itemize}

\vspace{1mm}

{\noindent{\bf\underline{Termination of the algorithm}.}}
In the case if the last coordinate of $\xi_{1,i}$ is zero
we do not proceed further and the algorithm terminates.

}

\begin{remark}
The algorithm was tested for quartic irrationalities, in all the examples the algorithm produced periodic output.
\end{remark}

\section{Dirichlet groups}\label{Dirichlet groups}

In order to understand better the reason of periodicity let us study maximal commutative subgroups of
$SL(d,\z)$. Such subgroups are called {\it Dirichlet groups}.

\subsection{Magic of integer commuting matrices}

Let us start with a simple exercise.

{\noindent{\bf Exercise 1}.} {\it Let}
$$
A=
\left(
\begin{array}{ccc}
2& 5& -1\\
3& 6& 1\\
4& 7& 1
\end{array}
\right).
$$
{\it Find an integer matrix with unit determinant commuting with $A$? }

\vspace{2mm}

The obvious solution to this exercise is the identity matrix, but let us disregard it.
Let us first peek the answer to this question.
The first matrix that we are able to detect is
$$
B=
\left(
\begin{array}{ccc}
88778750433916& 1881948516620816& -1642359549748757\\
-77918418013751& -849278651461089& 759124773173459\\
534000559063825& -721564227716990& 360094549931638
\end{array}
\right)
$$
The sizes of the elements of this matrix are rather impressive, are not they?
It is most likely that a brute force algorithm will reach the solution of this exercise only in the next millennium.

Even if one notices that the matrix $B$ is, in fact, a polynomial of $A$ with integer coefficients, namely
$$
B=-147205796095883A^2+1347947957556991A-399030223241821,
$$
the brute force search of the coefficients of such polynomial is still very long.
This is very much in the contrast to the complexity of the  input matrix $A$, each element of which  requires 4 bits only.
(It is also important to notice  that the coefficients of such polynomial are not necessarily integers,
they could be rational numbers instead.)

\vspace{2mm}

Let us discuss how to find the answer efficiently.

\subsection{Dirichlet groups}
First we give some necessary definitions.

\begin{definition}
Denote by $\Gamma(A)$ the set of all integer matrices commuting with $A$.

\vspace{1mm}

{\noindent
{\it $($i$)$} The {\it Dirichlet group} $\Xi(A)$ is the subset of invertible matrices in $\Gamma(A)$.
}

\vspace{1mm}

{\noindent
{\it $($ii$)$}
The {\it positive Dirichlet group}  $\Xi_+ (A)$ is the subset of $\Xi(A)$ that
consists of all matrices with positive real eigenvalues.
}
\end{definition}

\subsection{Dirichlet's unit theorem}\label{Dirichlet's unit theorem}
The first questions that one might ask before approaching Exercise~1 is
{\it whether such a matrix does exist}. Can it be that all unit determinant integer matrices commuting with $A$
are in fact powers of $A$? In terms of Dirichlet groups we ask wether the group $\Xi(A)$
is isomorphic to $\z$ or not.

The answer to this question is provided by the Dirichlet's unit theorem.
A precise formulation of the theorem is as follows.

\vspace{2mm}

{\noindent{ \bf Dirichlet's unit theorem.}}
{\it Let $K$ be a field of algebraic numbers of degree $n=s+2t$,
where $s$ is the number of real roots and $2t$ is the number of complex roots for the
minimal integer polynomial of any irrational elements of $K$.
Consider an arbitrary order $D$ in $K$.
Then $D$ contains units $\varepsilon_1,\ldots,\varepsilon_r$ for $r=s+t-1$
such that every unit
$\varepsilon$ in $D$ has a unique decomposition of the form
$$
\varepsilon=\xi\varepsilon_1^{a_1}\cdots\varepsilon_r^{a_r} \hbox{ },
$$
where $a_1,\ldots,a_r$ are integers and $\xi$ is a root of $1$ contained in $D$.}

\vspace{2mm}

So what is hidden behind Dirichlet's unit theorem? Rather than to go forward with all the formal definitions
involved in the formulation of this theorem
we prefer to reformulate this theorem simply in terms of matrices.
(For necessary definitions and the proof of the theorem we refer an interested reader, e.g.,
to the book~\cite{Borevich1966};
for the justification of the reformulation we refer, e.g., to Chapter 17 of~\cite{KarpenkovGCF2013}.
Also there is a lot of related material specifically on algebraic cubic and quartic fields in the book~\cite{Delone1964}.)

\vspace{2mm}

{\noindent{ \bf Dirichlet's unit theorem in the matrix form.}}
{\it
Let $A$ be an integer matrix $A$ whose characteristic
polynomial is irreducible over the field of rational numbers $\q$.
Let it has $s$ real and $2t$ complex eigenvalues.
Then there exists a finite Abelian group $G$ such that
$$
\Xi(A)=G\oplus \z^{s+t-1}.
$$
For the positive Dirichlet group we have:
$$
\Xi_+(A)=\z^{s+t-1}.
$$
}

\vspace{2mm}

\begin{example}
In three-dimensional case we have two possible situations.
\begin{itemize}
\item {\bf \underline{Complex case}:}
Let a characteristic polynomial has two complex roots.
In this case both Dirichlet and positive Dirichlet groups are isomorphic to $\z$.

\vspace{2mm}
\item {\bf \underline{Totally real case}:}
In case if all the roots of the characteristic polynomial are real numbers we have
$$
\Xi(A)=G\oplus \z^{2} \quad \hbox{and} \quad \Xi_+(A)=\z^{2}.
$$
\end{itemize}
\end{example}

\subsection{Several questions that we can answer}

The technique discussed in this paper gives a constructive approach the following questions.

\begin{question}\label{q1}
 Find any $\SL(3,\z)$-matrix commuting with a given integer matrix with irreducible characteristic polynomial over $\q$.
\end{question}

\begin{question}\label{q2}
Find an  $\SL(3,\z)$-matrix having a given cubic vector as an eigenvector.
\end{question}

\begin{question}\label{q3}
Let $M$ be any $\SL(3,\z)$-matrix whose characteristic polynomial is irreducible over $\q$.
Find an $\SL(3,\z)$-matrix commuting with $M$  that is not a power of $M$.
(Note that this question is interesting only in the totally real case as otherwise $\Xi(A)=\z$.)
\end{question}

In the next section we rewrite
Jacobi-Perron type algorithms in the matrix form and give the answers to these three questions.

\section{Jacobi-Perron type algorithms and Dirichlet groups}

\subsection{Jacobi-Perron algorithm in the matrix form}

Notice that Jacobi-Perron type algorithms can be formulated in terms of matrix multiplication form.
This concerns both the Jacobi-Perron algorithm,
the heuristic APD-algorithm,
the $\sin^2$-algorithm,
and many other similar algorithms (various types of such algorithms are collected in
the book of F~Schweiger~\cite{Schweiger2000}, see also Chapter~23.4 in~\cite{KarpenkovGCF2013}).
In case if an algorithm produces an eventually periodic output, the answers to Questions~1---3 above can be
obtained from the matrix form as we explain in the next subsection.

\subsection{Matrices with prescribed cubic eigenvectors}\label{Matrices with prescribed cubic eigenvectors}

Assume that the Jacobi-Perron algorithm is eventually periodic and its pre-period and period for a
given vector are respectively as follows:
$$
\left(
\left(
\begin{array}{c}
a_1\\
b_1\\
\end{array}
\right),
\ldots ,
\left(
\begin{array}{c}
a_n\\
b_n\\
\end{array}
\right)\right)
\quad \hbox{and} \quad
\left(
\left(
\begin{array}{c}
c_1\\
d_1\\
\end{array}
\right),
\ldots ,
\left(
\begin{array}{c}
c_m\\
d_m\\
\end{array}
\right)\right).
$$
Denote
$$
M_1=
\prod\limits_{i=1}^{n}
\left(
\begin{array}{ccc}
\mathbf{a_i}& 1& 0\\
1& 0& 0\\
\mathbf{b_i}& 0& 1
\end{array}
\right)
\quad \hbox{and} \quad
M_2=
\prod\limits_{j=1}^{n}
\left(
\begin{array}{ccc}
\mathbf{c_j}& 1& 0\\
1& 0& 0\\
\mathbf{d_j}& 0& 1
\end{array}
\right).
$$
and set
$$
M=M_1M_2(M_1)^{-1}.
$$
Then {\it $M$ is an $\SL(3,\z)$ matrix having original cubic vector as
the eigenvector whose  absolute value of the  eigenvalue is the greatest among all the absolute values of the eigenvalues of $M$}.

\begin{remark}
There are similar representations for both
the heuristic APD-algorithm and
the $\sin^2$-algorithm.
We omit them here, as they literally repeats the representation for the Jacobi-Perron algorithm
with obvious substitutions of matrices.
\end{remark}

\begin{example}
Let us discuss the vector
$$
(1,\xi,\xi^2+\xi)
$$
considered in  Example~\ref{ex-ok} above.
Recall that $\xi$ is a real root of the polynomial $x^3+2x^2+x+4$.
The Jacobi-Perron generates a periodic sequence with 6 steps of pre-period
and 2 steps of period:

\begin{center}
\begin{tabular}{|c||c|c|c|c|c|c||c|c||c|c||c|c||c|c|c|c|c|c|c|c|c|c|c}
\hline
&1 &2 &3 &4 &5 &6 &$2k+1$ &$2k+2$
\\
\hline
\hline
$\lfloor x/y\rfloor$ & -1 &1 &1 &1 &2 &6 &3 &7
\\
\hline
$\lfloor z/y\rfloor$& -2 &0 &0 &0 &2 &4 &1 &1
\\
\hline
\end{tabular}
\end{center}

First of all, we write matrices for a pre-period and a period:
\begin{align*}
M_1&=
\left(
\begin{array}{ccc}
\mathbf{-1}& 1& 0\\
1& 0& 0\\
\mathbf{-2}& 0& 1
\end{array}
\right)\cdot
\left(
\begin{array}{ccc}
\mathbf{1}& 1& 0\\
1& 0& 0\\
\mathbf{0}& 0& 1
\end{array}
\right)^3\cdot
\left(
\begin{array}{ccc}
\mathbf{2}& 1& 0\\
1& 0& 0\\
\mathbf{2}& 0& 1
\end{array}
\right)\cdot
\left(
\begin{array}{ccc}
\mathbf{6}& 1& 0\\
1& 0& 0\\
\mathbf{4}& 0& 1
\end{array}
\right)
\\
&=\left(
\begin{array}{ccc}
-22& -1& -3\\
51& 2& 7\\
-67& -3& -9
\end{array}
\right);
\\
M_2&=
\left(
\begin{array}{ccc}
\mathbf{3}& 1& 0\\
1& 0& 0\\
\mathbf{1}& 0& 1
\end{array}
\right)\cdot
\left(
\begin{array}{ccc}
\mathbf{7}& 1& 0\\
1& 0& 0\\
\mathbf{1}& 0& 1
\end{array}
\right)
=\left(
\begin{array}{ccc}
22& 1& 3\\
7& 0& 1\\
8& 0& 1
\end{array}
\right).
\end{align*}

Finally we get
$$
M=M_1M_2(M_1)^{-1}=
\left(
\begin{array}{ccc}
5& -4& 3\\
-12& 9& -7\\
16& -12& 9
\end{array}
\right).
$$
This concludes the computation of $M$.
\end{example}

\begin{example}
Let us consider once more the vector
$$
v=(1,\sqrt[3]{4},\sqrt[3]{16}) .
$$
As we have seen in Example~\ref{not-ok}
we are unable to get a periodic sequence generated by the Jacobi-Perron algorithm,
and hence we cannot find a requested matrix  using the Jacobi-Perron algorithm.

Let us use the heuristic APD-algorithm instead
(similarly multiplying the corresponding matrices for linear maps used in the algorithm).
Following the results of the continued fraction computations in  Example~\ref{not-ok-2} we find
that $v$ is an eigenvector of the matrix
$$
M=
\left(
\begin{array}{ccc}
5& 8& 12\\
3& 5& 8\\
2& 3& 5
\end{array}
\right).
$$
\end{example}

\subsection{Answers to Questions~\ref{q1}--\ref{q3}}
Finally we discus the answers to Questions~\ref{q1}--\ref{q3}.

\vspace{2mm}

{\noindent
{\bf Answer to Question~\ref{q1}.} }
Let us show how to construct an $\SL(3,\z)$-matrix commuting with a given integer matrix $A$ (assuming
that the characteristic polynomial of $A$ is irreducible over $\q$).
What we should do is to take a basis of eigenvectors of $A$ and run the heuristic APD-algorithm
for it (or $\sin^2$-algorithm in the totally real case).
The algorithm will generate a  (heuristically) eventually periodic sequence.
From its period and pre-period sequences one computes
the required matrix $M$ (as discussed in Subsection~\ref{Matrices with prescribed cubic eigenvectors}).
An example is shown in Exercise 1 above.

\vspace{2mm}

{\noindent
{\bf Answer to Question~\ref{q2}.} }
The second question is very similar to the first one.
Here we are requested to find an  $\SL(3,\z)$-matrix having a given cubic vector as an eigenvector.
Assume we are given by a cubic vector
$$
(1,\xi,q(\xi)),
$$
where $\xi$ is a solution of a certain cubic polynomial $p$.
Then we construct the other two conjugate vectors (following discussions of Remark~\ref{conj-remark}),
apply the heuristic APD-algorithm (or $\sin^2$-algorithm in the totally real case) and write the matrix from
the period and pre-period of the algorithm

\begin{remark}
Usually if we are given by a cubic vector, we are given by the corresponding polynomials $p$ and $q$.
However, this might be not the case, and we have an expression in a style of Kordano's formula instead.
In this case the polynomials $p$ and $q$ can be guessed from approximations
of $\xi^3$ and $q(\xi)$ and their approximate formulae in terms of approximations of $\xi^2$ and $\xi$ and $1$.
\end{remark}

\vspace{2mm}

{\noindent
{\bf Answer to Question~\ref{q3}.} }
Finally we write a matrix commuting with a given matrix $A\in SL(n\z)$ that is not a power of $A$.
As we have mentioned this question is valid only for the totally-real case
(in the complex case  the Dirichlet group of $A$
is isomorphic to $\z$).
Let $\xi$, $\nu$, and $\mu$ be eigenvectors of $A$.
First of all we construct $\SL(2,\z)$ matrices
$M_\xi$, $M_\nu$ and $M_\mu$ following the computations for the triples vectors
$$
(\xi, \nu, \mu), \quad
(\nu, \mu, \xi ), \quad \hbox{and} \quad
(\mu,\xi, \nu)
$$
respectively.
The main feature that we use further is that
{\it both the heuristic APD-algorithm, and $\sin^2$-algorithm
construct a matrix whose maximal absolute value of the eigenvalue
corresponds to the algebraic vector which is used in the iterations
in the first position}.

\vspace{2mm}

Notice that  maximal absolute values of the eigenvalues of $A^n$ correspond simultaneously to
the same eigenline for $n>0$; and for the same eigenline for $n<0$.
Hence there is one of the vectors $\xi$, $\nu$, and $\mu$ that is not on these two eigenlines.
Therefore,  the required matrix $M$ can be chosen from $M_\xi$, $M_\nu$ and $M_\mu$
by comparing the sizes of absolute values of eigenvalues.

\vspace{2mm}

Let us illustrate the answers to Questions~2 and~3 with the following example.

\begin{example}
Consider an irreducible cubic polynomial
$$
p(x)=2x^3 - 4x^2 - 7x - 2
$$
with three positive roots denoted by $\alpha$, $\beta$, and $\gamma$.
Our goal is to compute two independent (w.r.t. matrix multiplication)
$\SL(2,\z)$-matrices with eigenvectors
$$
\xi=(1,\alpha,\alpha^2), \quad
\nu=(1,\beta,\beta^2), \quad \hbox{and} \quad
\mu=(1,\gamma,\gamma^2).
$$
Direct computations using the heuristic APD-algorithm applied to triples
$$
(\xi, \nu, \mu), \quad
(\nu, \mu, \xi ), \quad \hbox{and} \quad
(\mu,\xi, \nu)
$$
result in the following three matrices:
$$
\begin{array}{c}
A=
\left(
\begin{array}{ccc}
55& 210& 176\\
176& 671& 562\\
562& 2143& 1795\\
\end{array}
\right);
\quad
B=
\left(
\begin{array}{ccc}
-497& -1122&  400\\
400& 903& -322\\
-322& -727& 259\\
\end{array}
\right);
\\
C=
\left(
\begin{array}{ccc}
185& 172& -72\\
-72& -67& 28\\
28& 26& -11\\
\end{array}
\right).
\end{array}
$$
Brute force search of the powers of matrices show that
$$
A^3B^5C^7=\id,
$$
where $\id$ is the identity matrix.
Since at least two of these matrices are independent, we have that any two of these matrices are linearly independent.
\end{example}

\begin{remark}
Clearly the triples of matrices $M_\xi$, $M_\nu$, and $M_\mu$ generate a finite
index sublattice in the positive Dirichlet groups. There exists a technique to find the basis using constructions of multidimensional
Klein continued fractions and observing combinatorics of their periods.
We do not discuss it in this paper and refer an interested reader to Chapter~20
of~\cite{KarpenkovGCF2013} (see also~\cite{Karpenkov2009}).
\end{remark}


\begin{remark}
The described in this section method works fine for $\SL(3,\z)$-matrices but it has some
limitations for $SL(d,\z)$-matrices with $d<3$.
For instance, in the totally real case of quartic case ($d=4$)
the corresponding positive Dirichlet group is three-dimensional.
A direct application of the current method potentially can output 4 matrices
spanning $\z^2$, and not $\z^3$ as is expected.
This should be a very rare case, {\it if possible at all}
(the last is unknown to the author).
\end{remark}

{\noindent
{\bf Acknowledgements.}
The author is grateful to H.~Servatius and A.Ustinov for useful comments and remarks.
}

\bibliographystyle{plain}
\bibliography{cfbiblio}

\vspace{.2cm}

\end{document}